\documentclass[12pt]{article}

\usepackage{amsmath, amsthm, amsfonts, amssymb,latexsym}
\usepackage{graphicx}
\usepackage[bookmarks=false]{hyperref}
\UseRawInputEncoding

\newtheorem{thm}{Theorem}[section]
\newtheorem{cor}[thm]{Corollary}
\newtheorem{lem}[thm]{Lemma}

\theoremstyle{definition}
\newtheorem{defn}[thm]{Definition}

\begin{document}

\def\prof{{\sc Proof.\ \ }}
\def\sect#1{\begin{center}\section{#1}\end{center}} 
\def\R#1{{\bf R}^{#1}}
\def\I#1#2{\int\limits_{#1}#2} 
\def\p#1{#1^{\prime}} 
\def\l#1#2#3{\lim_{#1 \rightarrow #2}#3}
\def\ld#1#2#3{\liminf_{#1 \rightarrow #2}#3} 
\def\lu#1#2#3{\limsup_{#1 \rightarrow #2}#3}
\def\E#1#2{{\bf E}^{#1}\left(#2\right)} 
\def\P#1#2{{\bf P}^{#1}\left(#2\right)} 
\def\Pt#1#2#3{{\bf P}(#1,#2,#3)} 
\def\pb#1{{\bf P}\left(#1\right)}
\def\eb#1{{\bf E}\left(#1\right)}
\def\Di#1{\,Dim\,#1}
\def\di#1{\,dim\,#1}
\def\fr#1#2{\frac{#1}{#2}}
\def\beq{\begin{equation}}
\def\eeq{\end{equation}}
\def\bea{\begin{eqnarray}}
\def\bean{\begin{eqnarray*}}
\def\eean{\end{eqnarray*}}
\def\eea{\end{eqnarray}}
\def\heq#1#2#3{\hbox to \hsize{\hskip #1 $#2$ \hss (#3)}}
\def\df#1#2{\frac{\displaystyle #1}{\displaystyle #2}}

\def\bdes{\begin{description}}
\def\ndes{\end{description}}

\newcommand{\bh}{{\bf h}}
\newcommand{\hf}{{\bf f}}
\newcommand{\he}{{\bf e}}
\newcommand{\hL}{{\bf L}}
\newcommand{\hg}{{\bf g}}
\newcommand{\hG}{{\bf G}}
\newcommand{\hM}{{\bf M}}

\newcommand{\bb}{\mathbb{b}}
\newcommand{\ww}{\mathbb{W}}
\newcommand{\hh}{\mathbb{H}}
\newcommand{\dd}{\mathbb{D}}
\newcommand{\cc}{\mathbb{C}}
\newcommand{\ee}{\mathbb{E}}
\newcommand{\zz}{\mathbb{Z}}
\newcommand{\nn}{\mathbb{N}}
\newcommand{\pp}{\mathbb{P}}
\newcommand{\qq}{\mathbb{Q}}
\newcommand{\ttt}{\mathbb{T}}

\def\rr{\Bbb R}
\def\a{\alpha}
\def\b{\beta}
\def\g{\gamma}
\def\s{\sigma}
\def\ep{\epsilon}
\def\d{\delta}
\def\D{\Delta}
\def\o{\omega}
\def\O{\Omega}
\def\n={\not=}
\def\u>{\wedge}
\def\d>{\vee}
\def\.{\bullet}
\def\l{\lambda}
\def\L{\Lambda}
\def\r{\rho}
\def\vf{\varphi}
\def\f{\phi}
\def\t{\tau}
\def\z{\zeta}
\def\na{\nabla}

\def\B{{I\!\!B}}
\def\E{{I\!\!E}}
\def\N{{I\!\!N}}
\def\P{{I\!\!P}}
\def\Q{{I\!\!\!Q}}
\def\R{{I\!\!R}}

\def\AA{\mathcal A}
\def\BB{\mathcal B}
\def\CC{\mathcal C}
\def\DD{\mathcal D}
\def\FF{\mathcal F}
\def\EE{\mathcal E}
\def\GG{\mathcal G}
\def\JJ{\mathcal J}
\def\LL{\mathcal L}
\def\NN{\mathcal N}
\def\PP{\mathcal P}
\def\SS{\mathcal S}

\def\Fht{\hbox{${\cal F}_{t}$}}
\def\Ght{\hbox{${\cal G}_{t}$}}

\def\U{\bigcup}
\def\Uu{\bigcap}
\def\Au{\forall}
\def\Eu{\exists}
\def\8u{\infty}
\def\0{\emptyset}
\def\mp{\longmapsto}
\def\rt{\rightarrow}
\def\up{\uparrow}
\def\dn{\downarrow}
\def\sub{\subset}
\def\vep{\varepsilon}
\def\ln{\langle}
\def\rn{\rangle}

\def\<<{\langle\!\langle}
\def\>>{\rangle\!\rangle}
\def\3|{|\!|\!|}
\def\ess{\text{\rm{ess}}}
\def\beg{\begin}

\def\beqt{\begin{equation}}
\def\neqt{\end{equation}}
\def\beq{\begin{equation}}
\def\neq{\end{equation}}

\def\bdes{\begin{description}}
\def\ndes{\end{description}}

\def\Ric{\text{\rm{Ric}}}
\def\Hess{\text{\rm{Hess}}}
\def\i{\text{\rm{i}}}
\def\ii{\text{\rm{ii}}}
\def\iii{\text{\rm{iii}}}
\def\iv{\text{\rm{iv}}}
\def\v{\text{\rm{v}}}
\def\vi{\text{\rm{vi}}}
\def\vii{\text{\rm{vii}}}
\def\viii{\text{\rm{viii}}}
\def\e{\text{\rm{e}}}

\def\Ra{\Rightarrow}
\def\Lra{\Leftrightarrow}
\def\rto{\longrightarrow}
\def\la{\leftarrow}
\def\ra{\rightarrow}
\def\ua{\uparrow}
\def\da{\downarrow}


\title{The existence of a strongly polynomial time simplex algorithm for linear programs
\thanks{Supported by the National Natural Science Foundation of China (11871118,12271061).}}

\author{Zi-zong Yan, Xiang-jun Li and Jinhai Guo \thanks{School of Information and Mathematics,
Yangtze University, Jingzhou, Hubei,
China(zzyan@yangtzeu.edu.cn, franklxj001@163.com and xin3fei@21cn.com).}}
\date{}
\maketitle

\begin{abstract} It is well known that the most challenging question in optimization and discrete geometry is whether there is a strongly polynomial time simplex algorithm for linear programs (LPs). This paper gives a positive answer to this question by using the parameter analysis technique presented by us (http://arxiv.org/abs/2006.08104). We show that there is a simplex algorithm whose number of pivoting steps does not exceed the number of variables of a LP problem.

\textbf{Keywords:} Linear programming, parametric linear programming, complement slackness property, projection polyhedron, simplex algorithm, pivot rule, polynomial complexity \\

\textbf{AMS subject classifications. } 90C05, 90C49, 90C60, 68Q25, 68W40 
\end{abstract}

\section{Introduction}
Linear programming (LP) is the problem of minimizing a linear objective function over a polyhedron $P\subset\mathbb{R}^n$ given by a system of $m$ equalities and $n$ nonnegative variables. As Chakraborty et al. \cite{CTA20} recently put it,

{\it "The story of LP is one with all the elements of a grand historical drama. The original idea of testing if a polyhedron is non-empty by using a variable elimination to project down one dimension at a time until a tautology emerges dates back to a paper by Fourier \cite{Fou26} in 1823. This gets re-invented in the 1930s by Motzkin \cite{Mot36}. The real interest in LP happens during World War II when mathematicians ponder best ways of utilising resources at a time when they are constrained". }

\noindent As one of the fundamental problems in optimization and discrete geometry, LP has been a very successful undertaking in the field of polyhedral combinatorics in the last six decades. Part of this success relies on a very rich interplay between geometric and algebraic properties of the boundary of a polyhedron and corresponding combinatorial structures of the problem it encodes.

The simplex algorithm, invented by Dantzig \cite{Dan47} solves LPs by using pivot rules and procedures an optimal solution. Subsequently developed various simplex algorithms are to iteratively improve the current feasible solutions by moving from one vertex of the polyhedron to an adjoint one according to some pivot rule, until no more improvement is possible and optimality is proven. All edges that connect these adjoint vertices of the polyhedron form a pivot path. The complexity of the simplex algorithm is then determined by the length of the path - the number of pivot steps.

It is a theoretically interesting and practically relevant question to recognize how a pivot rule \cite{Sha87} quickly leads to the optimal vertex. Connected to this are the diameter problem and the algorithm problem, for example, see \cite{GHZ98}. Closely to the Hirsch conjecture and its variants, the diameter problem is whether there is a short path to the optimal vertices. Unfortunately Santos' recent counter-example \cite{San12} disproves this conjecture. On the subject for more information see the papers \cite{Bor87,GK93,Mef87,Mir12} and the comments \cite{Eis13,HU13,Lam13,San13,Ter13}. The algorithm problem is whether there is a very strongly polynomial time algorithm for LP. There has been recent interest in finding an algorithm like this for some  special cases, such as the deterministic Markov decision processes \cite{PY13}, the generalized circulation problem \cite{GJL02}, the maximum flow problem \cite{AW98,GC97}, the minimum cost network flow problem \cite{AW97}, the Arrow-Debreu Market Equilibrium \cite{Jai07} and so on. However, up to now, no pivoting rule yields such a algorithm depending only on the polynomiality of the number of constraints and the number of variables in general, for example, see \cite{Ter13}. In contrast to the excellent practical performance of the simplex algorithm, the worst-case time complexity of each analyzed pivot rule in the field of LP is known to grow exponentially, for example, see \cite{AC78,Gol83,Gol79,KL72,MPZ01,Mur76,Mur80,Pap89,Roo90}. Other exponential example was presented by Fukuda and Namiki \cite{FN94} for linear complementarity problems.

The simplex algorithm belongs to the '10 algorithms with the greatest influence on the development and practice of science and engineering in the 20th century', see \cite{DS00}. It has performed sufficiently well in practice, but theoretically the complexity of a pivot rule has so far been quite a mystery. To explain the large gap between practical experience and the disappointing worstcase, the tools of average case analysis and smoothed analysis have been devised, and to conquer the worst case bounds, research has turned to randomized pivot rules. For the average case analysis, a polynomial upper bounded was achieved \cite{Bor87}. In contrast, so far, smoothed analysis has only been done for the shadow-vertex pivot rule introduced by Gaas and Saaty \cite{GS55}. Under reasonable probabilistic assumptions its expected number of pivot steps is polynomial \cite{Spi01}. It is worth noting that none of the existing results exclude the possibility of (randomized) pivot rules being the desired (expected) polynomial-time pivot rules. For more information on the randomized pivot rules reference the papers \cite{EV17,GV07,GHZ98,Kal92,MSW,Mur80,Ver09} and etc.

Two important advances have been made in polynomial time solvability for the ellipsoid method developed by Khachain \cite{Kha80} and the interior-point method initiated by Karmarkar \cite{Kar84} since the 1970s. However, the run time¡¯s complexity of such two algorithms is only qualified as weak polynomial. Actually, constructing a strongly polynomial time pivot rule is the most challenging open question in optimization and discrete geometry \cite{GK93,Mef87,Sma00,Tod00}.

The survey by Terlaky and Zhang \cite{TZ93} contributes the various pivot rules of the simplex algorithm and its variants, in which they categorized pivot rules into three types. The first alternative, which is called combinatorial pivot rule, is to take care of the sign of the variables, including either costs or profits. The algorithms of this type known to the authors are that of Bland \cite{Bla77}, Folkman and Lawrence \cite{FL78}, Fukuda \cite{Fuk82} and etc. The second alternative, which is called parametric pivot rule, is closely related to parametric programming, more precisely to the shadow vertex algorithm \cite{Boy92,Gol83,Mur80} and to Dantig's self-dual parametric simplex algorithm \cite{Dan63}. The algorithms can be interpreted as Lemke's algorithm \cite{Lem65} for the corresponding linear complementarity problem \cite{Lus87}. Algorithms of this third type, which are close connections to certain interior point methods (see Karmarker \cite{Kar84}), allow the iterative points to go inside the polytope. It is believed to be able to use more global information and therefore to avoid myopiness of the simplex method, for example, see Roos \cite{Roo86}, Todd \cite{Tod90} and Tamura et al. \cite{TTF88}.

Our main goal of this paper is to give a positive answer to the above open question by the use of the parametric analysis technique that we recently proposed in \cite{YLG20,YLG22}. We show that there exists a simplex algorithm whose number of pivoting steps does not exceed the number of variables of a LP problem.

The organization of the rest of the paper is as follows. In Section 2, we recall the strong duality theorem of LP. In Section 3, we review the main results of the paper \cite{YLG20,YLG22}: we define the set-valued mappings between two projection polyhedrons for parametric LP problems and establish the relationship between perturbing the objective function data (OFD) and perturbing the right-hand side (RHS) for a LP problem without using dual. We then investigate the projected behavior of the set-valued mappings in Section 4. As a application, we present the existence of a strongly polynomial time pivot rule for a LP   in the final section.

\section{Preliminaries} \label{sect21}
Consider the following pair of LPs in the standard forms:
\begin{equation} \label{primaljihe10p} \begin{array}{ll} \min\limits_x & \langle c,x \rangle \\ s. t. & Ax=b,\\ & x \geq 0 \end{array} \end{equation} and
\begin{equation} \label{primaljihe10d} \begin{array}{ll} \max\limits_w & b^T w \\ s. t. & A^Tw\leq c, \end{array} \end{equation} where $c\in\mathbb{R}^n, b\in\mathbb{R}^m$ and $A\in\mathbb{R}^{m\times n}$ are given. As usual, the first LP (\ref{primaljihe10p}) is called the primal problem and the second problem (\ref{primaljihe10d}) is called the dual problem, and the vectors $c$ and $b$ are called the cost and the profit vectors, respectively. By $P=\{x\in\mathbb{R}^n| Ax=b, x\geq 0\}$ and $D=\{w \in\mathbb{R}^m | A^Tw\leq c\}$, we denote the feasible sets of the primal and dual problems, respectively.

The strong duality theorem provides a sufficient and necessary condition for optimality, for example, see \cite{Dan47,GT56,Mur76,Sch86}.

\thm \label{weakd2} {\bf (Strong duality)} If the primal-dual problem pair (\ref{primaljihe10p}) and (\ref{primaljihe10d}) are feasible, then the two problems are solvable and share the same objective value. \upshape

If the primal and dual programs have optimal solutions and the duality gap is zero, then the Karush-Kuhn-Tucker (KKT) property for the primal-dual LP problem (\ref{primaljihe10p}) and (\ref{primaljihe10d}) pair is
\begin{subequations}
\begin{align} \label{jbkkt1} & Ax = b, \quad x\geq 0, \\ \label{jbkkt2} & A^Tw\leq c, \\ \label{jbkkt3} & \langle x,c-A^Tw \rangle=0,
\end{align}
\end{subequations} in which the last equality is called the complement slackness property. Conversely, if $(x^*,w^*)\in\mathbb{R}^n\times \mathbb{R}^m$ is a pair of solutions of the system (\ref{jbkkt1})-(\ref{jbkkt3}), then $(x^*,w^*)$ is a pair of optimal solutions of the pair of primal-dual problems (\ref{primaljihe10p}) and (\ref{primaljihe10d}).

The feasible region of a LP problem is a polyhedron. In particular, a polyhedron is called a polytope if it is bounded. A nontrivial face $F$ of a polyhedron is the intersection of the polyhedron with a supporting hyperplane, in which $F$ itself is a polyhedron of some lower dimension. If the dimension of $F$ is $k$ we call $F$ a $k$-face of the polyhedron. The empty set and the polyhedron itself are regarded as trivial faces. $0$-faces of the polyhedron are called vertices and $1$-faces are called edges. Two different vertices $x^1$ and $x^2$ are neighbors if $x^1$ and $x^2$ are the endpoints of an edge of the polyhedron. That is, the linear segment $[x^1,x^2]=\{\lambda x^1+(1-\lambda)x^2|0\leq \lambda\leq 1\}$ is an edge of $P$. For material on convex polyhedron and for many references see Ziegler's book \cite{Zie95}.

$x\in\mathbb{R}^n$ is called a basic feasible solution of the primal problem (\ref{primaljihe10p}) if $x$ is a vertex of $P$, and $y=c-A^Tw\in\mathbb{R}^n$ is called a basic feasible solution of the dual problem (\ref{primaljihe10d}) if $w$ is a vertex of $D$.

\section{Parametric KKT conditions}
Describing a pivot path of the simplex algorithm is a difficult task. The ingredient in our construction is to investigate the relationship between two projection polyhedrons associated with a pair of almost primal-dual problems. Once the relationship is available, either the OFD or the RHS perturbations of the primal and the dual problems are closely linked. The concepts and results in this section come all from the paper \cite{YLG20,YLG22}.

Consider the following two parametric LP problems with two independent parametric vectors $u,v \in \mathbb{R}^r$
\begin{equation} \label{primaljihe1} \begin{array}{ll} \min\limits_x & \langle c+M^Tu,x \rangle \\ s. t. & Ax=b,\\ & x \geq 0, \end{array} \end{equation} and
\begin{equation} \label{primaljihe2} \begin{array}{ll} \max\limits_{w,s} & b^Tw+(M^Td+v)^Ts \\ s. t. & A^Tw+M^Ts\leq c, \end{array} \end{equation}
where $A\in\mathbb{R}^{m\times n}, M \in\mathbb{R}^{r\times n}$, $c,d\in \mathbb{R}^n$ and $b\in \mathbb{R}^m$ are given.
The following technical claims are assumed to hold throughout this paper.

{\bf Assumption 1}. The problems (\ref{primaljihe1}) and (\ref{primaljihe2}) are feasible.

{\bf Assumption 2}. $MM^T=I_r$, a $r\times r$ unit matrix.

{\bf Assumption 3}. The range spaces $R(A^T)$ and $R(M^T)$ are orthogonal, where $R(A^T)$ denotes the range space of $A^T$.

It should be noted that column vectors of $M^T$ denote perturbation directions of the objective function of the problem (\ref{primaljihe1}). On the other hand, authors \cite{YLG20} assumed that $d$ and $(0,0)$ are feasible for the problems (\ref{primaljihe1}) and (\ref{primaljihe2}), respectively. Although this assumption is unnecessary for this paper, for convenience, we still adopt the same notation of the paper \cite{YLG20} to preserve the vector $d$.

For the parametric LPs (\ref{primaljihe1}) and (\ref{primaljihe2}), their Lagrangian dual problems can be expressed as follow
 \begin{equation} \label{primaljihe1dual} \begin{array}{ll} \max\limits_w & b^Tw \\ s. t. & A^Tw\leq c+M^Tu \end{array} \end{equation}
and \begin{equation} \label{primaljihe2dual} \begin{array}{ll} \min\limits_x & \langle c, x \rangle \\ s. t. & Ax=b, \\ & Mx=Md+v,\\ & x \geq 0, \end{array} \end{equation} respectively. Clearly, the perturbations in the problems (\ref{primaljihe1dual}) and (\ref{primaljihe2dual}) occur in the RHS and not in the OFD.


The problems (\ref{primaljihe1}) and (\ref{primaljihe2dual}) represent two different perturbations of the problem (\ref{primaljihe10p}). The perturbations in the two problems occur in the OFD and in the RHS, respectively. To characterize the relationship, let us define one set-valued mapping as follows
\[ \Phi(u)= \{ M(x^*(u)-d)| x^*(u) \ \operatorname{is \ an \ optimal \ solution\ of}\ (\ref{primaljihe1}) \}. \] Analogously, we may define another set-valued mapping as follows
\[ \Psi(v)= \{ -s^*(v)| (w^*(v),s^*(v)) \ \operatorname{is \ a \ pair \ of \ optimal \ solution\ of}\ (\ref{primaljihe2}) \}. \]
Here the value of the mapping $\Phi(u)$ (or $\Psi(v)$) could be a set if the optimal solution is not unique corresponding to the parameter $u$ (or $v$). The similar idea was used by Borrelli et al. \cite{BBM01}.

Clearly, the domains of the set-valued mappings $\Psi$ and $\Phi$ are contained in the following two sets
\begin{eqnarray*} \varTheta_P &=& \{v\in\mathbb{R}^r| Ax=b, Mx=Md+v, x\geq 0\}, \\
\varTheta_D &=& \{u\in\mathbb{R}^r| \exists w\in\mathbb{R}^m \ \operatorname{ s.t. }\ A^Tw\leq c+M^Tu\}. \end{eqnarray*}
Obviously, the set $\varTheta_D$ is a projection of the polyhedron $\{(u,w)\in\mathbb{R}^{r+m}| c+M^Tu-A^Tw\geq 0\}$.
The following result shows that the set $\varTheta_P$ is a projections of another polyhedron.

\cor\label{xinzd1} Let $B\in\mathbb{R}^{n\times l}$ such that $R(A^T),R(M^T),R(B^T)$ span the whole space $\mathbb{R}^n$, and they are orthogonal with each other. Then there is a vector $t\in\mathbb{R}^l$ such that the primal slack vector \[x=d+M^Tv+B^Tt\] is feasible for the problem (\ref{primaljihe2dual}) if and only if $v\in \varTheta_P$. \upshape

\prof References the paper \cite[Corollary 2.4]{YLG20}. $\square$

It's worth noting that in this corollary, the choice of the matrix $B$ is not unique.

The set-valued mappings $\Phi$ closely links two different kinds of perturbations of (\ref{primaljihe10p}). Geometrical, the OFD perturbation means that the hyperplane $H=\{x\in \mathbb{R}^n| -\langle c+M^Tu,x\rangle=-\langle c+M^Tu,x^*(u)\rangle\}$ passing through $x^*(u)$ supports the feasible region $P$, and the RHS perturbation means that the affine set $C=\{x\in\mathbb{R}^n| Mx=Md+v\}$ cuts the region $P$, in which the affine set $C$ passes through $x^*(u)$ if $v\in\Phi(u)$. A similar geometric explanation applies the set-valued mappings $\Psi$.

\thm\label{maithm1} (1) For every $u\in \varTheta_D$, $\Phi(u)$ is well defined such that $\Phi(u) \subset \varTheta_P; $

(2) For every $v\in \varTheta_P$, $\Psi(v)$ is well defined such that $ \Psi(v) \subset \varTheta_D. $ \upshape

\prof We only prove the first result. By Theorem \ref{weakd2}, for every $u\in \varTheta_D$, the problem (\ref{primaljihe1}) is solvable that $\Phi(u)$ is well defined. If $x^*(u)$ is an optimal solution of the problem (\ref{primaljihe1}) and if $v\in\Phi(u)$, then $x^*(u)$ is also an optimal solution of the problem (\ref{primaljihe2dual}). Then by Corollary \ref{xinzd1}, one has $\Phi(u) \subset \varTheta_P$. $\square$

The set-valued mappings $\Phi$ and $\Psi$ play an important role in the development of strong duality theory and in the sensitivity analysis of conic linear optimization. An important and useful result is given in the paper \cite{YLG20,YLG22}.

\thm \label{interior2yy} Let $u\in \varTheta_D$ and $v\in\varTheta_P$ be arbitrary. Then there are a triple of vectors $(\bar{x},\bar{w},\bar{s})\in\mathbb{R}^{n+m+r}$ with $\bar{s}=-u$ such that the multiparametric KKT (mpKKT) property holds \begin{subequations}
\begin{align}\label{kktcond1} &A \bar{x} = b,\quad M \bar{x} =Md+v, \quad \bar{x}\geq 0, \\ \label{kktcond2} &
A^T\bar{w}+M^T\bar{s} \leq c, \\ \label{kktcond3} & \langle \bar{x}, c-A^T\bar{w}-M^T\bar{s} \rangle=0,
\end{align}
\end{subequations} if and only if both $v\in\Phi(u)$ and $u\in \Psi(v)$ hold. Furthermore, $(\bar{x},\bar{w})$ is a pair of optimal solutions of the problems (\ref{primaljihe2dual}) and (\ref{primaljihe1dual}) and $(\bar{x};\bar{w},\bar{s})$ is a pair of optimal solutions of the problems (\ref{primaljihe1}) and (\ref{primaljihe2}). \upshape

\prof References the paper \cite[Theorem 3.10]{YLG20}. $\square$

The above mpKKT property (\ref{kktcond1})- (\ref{kktcond3}) extends the classical KKT property (\ref{jbkkt1})-(\ref{jbkkt3}). More detail of remarks can be found in the paper \cite{YLG20,YLG22}. For this reason, either the problems (\ref{primaljihe2dual}) and (\ref{primaljihe1dual}) or (\ref{primaljihe1}) and (\ref{primaljihe2}) are called a pair of primal and dual LPs. Corresponding, $\varTheta_P$ and $\varTheta_D$ are called a pair of primal and dual projection polyhedra.

\thm \label{xinzd4} One has \begin{subequations} \begin{align} \label{jihexd1} \bigcup\limits_{u\in \varTheta_D} \Phi(u)= \varTheta_P,\\ \label{jihexd2} \bigcup\limits_{v\in \varTheta_P} \Psi(v) = \varTheta_D.
\end{align}
\end{subequations} \upshape

\prof By Theorem \ref{maithm1}, we have
\begin{equation} \label{nbhgx1} \bigcup\limits_{u\in \varTheta_D} \Phi(u)\subset \varTheta_P. \end{equation} Following we prove that the inverse inclusion relation is true.

Note that the feasible region of of the problem (\ref{primaljihe2}) does not depend on the given parametric vector $v$. Let us number all feasible bases of the problem (\ref{primaljihe2}) and let $\mathcal{I}$ denotes an index set
\[ \mathcal{I} =\{i| \exists (w^i,s^i) \operatorname{\ s.t.\ } y^i =c-A^Tw^i-M^Ts^i \ \operatorname{is \ a \ feasible\ basis\ of\ the \ dual\ problem}\ (\ref{primaljihe2})\}. \]
Then for every $i\in\mathcal{I}$, $y^i$ is a dual feasible basis and $\hat{x}^i(v)$ is the corresponding dual cost (a primal basis, perhaps $\hat{x}^i(v)$ is primal infeasible for some $v\in \mathbb{R}^r$) such that the complementary condition holds, i. e., $\langle \hat{x}^i(v),y^i\rangle=0$. Defined the set
\[\varTheta_i=\{v\in\mathbb{R}^r| \hat{x}^i(v)\geq 0 \}. \] Clearly, one has $\varTheta_i\subset \varTheta_P$ and $\bigcup\limits_{i\in \mathcal{I}}\varTheta_i=\varTheta_P$. If $\varTheta_i$ is not empty, then for every $v\in\varTheta_i$, $\hat{x}^i(v)$ is feasible and is also optimal for the problem (\ref{primaljihe2dual}). Therefore, the set $F_i=\{\hat{x}^i(v)|v\in\varTheta_i\}$ is a nontrivial face of the feasible region of the primal problem (\ref{primaljihe2dual}); and $y^i$ is also optimal for the problem (\ref{primaljihe1dual}) corresponding to the parameter $u=u^i=-s^i$. Furthermore, by Theorem \ref{interior2yy}, one has $\Phi(u^i)=\varTheta_i$ such that the inverse of the inclusion relation (\ref{nbhgx1}) holds. $\square$

\cor \label{adjhi} Let $u\in \varTheta_D$. Then $v\in\Phi(u)$ if and only if $u\in\Psi(v)$. \upshape

\prof Assume that $v\in\Phi(u)$. By Theorem \ref{xinzd4}, one has $v\in \varTheta_P$. Then the problems (\ref{primaljihe1}) and (\ref{primaljihe1dual}) are solvable, and the problem (\ref{primaljihe2dual}) is also solvable (see \cite[Corollary 3. 3]{YLG20}). Furthermore, the problem (\ref{primaljihe2}) is solvable and $u\in\Psi(v)$. And vice versa. $\square$

\section{The projection transformation}
Consider a pair of primal-dual parametric LPs
\begin{equation} \label{primaljihe10c} \begin{array}{ll} \min\limits_x & \langle c+M^TSu,x \rangle \\ s. t. & Ax=b,\\ & x \geq 0 \end{array} \end{equation} and
\begin{equation} \label{primaljihe10cd} \begin{array}{ll} \max\limits_{w,s} & b^Tw+(Md+Sv)^Ts \\ s. t. & A^Tw+M^Ts\leq c, \end{array} \end{equation}
 where $S\in\mathbb{R}^{r\times r}$ is a symmetric projection matrix, i. e., $S^T=S=S^2$. Such a pair of primal-dual LPs denote the projection of the pair of parametric primal-dual LPs (\ref{primaljihe1}) and (\ref{primaljihe2}).

Similarly, the following mpCLO problems
\begin{equation} \label{primaljihe40c} \begin{array}{ll} \min\limits_x & \langle c, x \rangle \\ s. t. & Ax=b,\\ & SMx=SMd+ Sv, \\ & x \geq 0 \end{array} \end{equation}
and \begin{equation} \label{primaljihe20c} \begin{array}{ll} \max\limits_w & b^Tw \\ s. t. & A^Tw\leq c+M^TSu \end{array} \end{equation} are the Lagrangian duals of problems (\ref{primaljihe10c}) and (\ref{primaljihe10cd}), respectively. Then the corresponding mpKKT property is as follows
\begin{subequations}
\begin{align}\label{skktcond1} &A x =b,\quad Mx= Md + Sv,\quad x\geq0 , \\ \label{skktcond2} &
A^Tw\leq c+M^TSu,\\ \label{skktcond3}& \langle x, c+M^TSu-A^Tw \rangle=0.
\end{align}
\end{subequations}

\thm\label{mfzjia1} Let $S\in\mathbb{R}^{r\times r}$ be a symmetric projection matrix.

(1) For every $u\in \varTheta_D$, there exists a vector $v\in S\varTheta_P=\{Sv|v\in \varTheta_P\}$ such that $Sv\in\Phi(Su)$ and $Su\in\Psi(Sv)$;

(2) For every $v\in \varTheta_P$, there exists a vector $u\in S\varTheta_D=\{Su|u\in \varTheta_D\}$ such that $Su\in\Psi(Sv)$ and $Sv\in\Phi(Su)$. \upshape

\prof Clearly, if $v\in S\varTheta_P$ or $u\in S\varTheta_D$, then $Sv=v$ or $Su=u$. Then the result follows from Theorem \ref{maithm1} and Corollary \ref{adjhi}. $\square$

The more information about this sunject references the paper \cite{YLG22}.

\section{Polynomial complexity}
Prior to the study of sensitivity analysis in recent years, the actual invariancy region plays an important role in the development of parametric LP. Adler and Monteiro \cite{AM92} first investigated the sensitivity analysis of LPs by using the optimal partition approach, in which they identify the range of parameters where the optimal partition remains invariant. Other treatments of the sensitivity analysis for LPs based on the same approach was given by Jansen et al. \cite{JR03}, Greenberg \cite{Gre94}, and Roos et al. \cite{RTV05}, Ghaffari-Hadigheh et al. \cite{GG08}, Berkelaar et al. \cite{Ber97}, Dehghan et al. \cite{Deh07}, Hlad$\acute{i}$k \cite{Ha10} and etc.

Recently we \cite{YLG20} developed a novel optimal partition approach for conic linear optimization.

\defn \label{djkyi1} Let $\mathcal{V}$ be a simply connected subset of $\varTheta_P$. Then $\mathcal{V}$ is called an invariancy set if for all $v^1,v^2\in\mathcal{V}$, one of the following statements holds:

 (1) $\Psi(v^1)=\Psi(v^2)$ if $\mathcal{V}$ is not a singleton set;

 (2) $\Psi(v^1)$ is not a singleton set if $\mathcal{V}$ is equal to $\{v^1\}$.

In the approach one can identify the range of parameters where the optimal partition invariant. In this definition, the first claim means that the dual optimal objective value remains unchange, whereas the second claim means that the primal optimal objective value remains unchange.

The notation $\operatorname{dim}(\mathcal{V})$ denotes the dimension of the affine hull of the set $\mathcal{V}$.

\defn Let $\mathcal{V}$ be an invariancy set of $\varTheta_P$. Then $\mathcal{V}$ is called a transition face if $\operatorname{dim}(\mathcal{V})<r$. In particular, if $\operatorname{dim}(\mathcal{V})=0$, then $\mathcal{V}$ is called a transition point; and if $\operatorname{dim}(\mathcal{V})=1$, then $\mathcal{V}$ is called a transition line, and etc. \upshape

For the dual projection polyhedron $\varTheta_D$, the definitions of the invariant set and the transition face is similar.


\thm \label{zyjg1} Any two different invariancy regions of a projection polyhedron do not intersect. \upshape

\prof This result follows immediately from Definition \ref{djkyi1}. $\square$

Following we present several auxiliary results for the single-parameter LPs, in which a nontrivial invariancy set is a open interval; moreover, the endpoints of invariancy interval are the transition points (see also \cite{RTV05}).

\lem \label{AMyiny2} Let $r=1$.

(1) If $u=-\bar{s}\in\varTheta_D$ is given, then $\Phi(u)=[\underline{v},\overline{v}]$ can be identified by solving the following two auxiliary LP problems:
\begin{eqnarray*} \underline{v} &=& \min\{v| \operatorname{\exists \ (\bar{x};\bar{w})\ s. t. \ the\ parametric\ KKT\ property } \ (\ref{kktcond1})-(\ref{kktcond3})\ \operatorname{holds}\},\\ \overline{v}&=& \max \{v| \operatorname{\exists \ (\bar{x};\bar{w})\ s. t. \ the\ parametric\ KKT\ property } \ (\ref{kktcond1})-(\ref{kktcond3})\ \operatorname{holds}\}. \end{eqnarray*}

(2) If $v\in\varTheta_P$ is given, then $\Psi(v)=[\underline{u},\overline{u}]$ can be identified by solving the following two auxiliary LP problems:
\begin{eqnarray*} \underline{u} &=& \min\{-\bar{s}| \operatorname{\exists \ (\bar{x};\bar{w},\bar{s})\ s. t. \ the\ parametric\ KKT\ property }\ (\ref{kktcond1})-(\ref{kktcond3})\ \operatorname{holds}\},\\ \overline{u}&=& \max \{-\bar{s} | \operatorname{\exists \ (\bar{x};\bar{w},\bar{s})\ s. t. \ the\ parametric\ KKT\ property } \ (\ref{kktcond1})-(\ref{kktcond3})\ \operatorname{holds}\}. \end{eqnarray*} \upshape
\indent \prof The result follows from Theorem \ref{interior2yy}. $\square$

The open interval $(\underline{v},\overline{v})$ is an invariancy interval, in which the finite endpoints of the invariancy interval are transition points. Of course, either $\underline{v}=-\infty$ or $\overline{v}=+\infty$ is allowed. The similar results apply the set-valued mapping $\Psi$. Two direct consequences of Lemma \ref{AMyiny2} are as follow.

\lem \label{AMyiny3} Let $r=1$.

(1) If $\Phi(u)=[\underline{v},\overline{v}]$ is an interval, where $\underline{v}\ne\overline{v}$, then all optimal solutions $\hat{x}^*(v)$ of the problem (\ref{primaljihe2dual}) for $v\in [\underline{v},\overline{v}]$ forms an edge of the primal feasible region $P$, in which $\hat{x}^*(\underline{v})$ and $\hat{x}^*(\overline{v})$ are two adjoint vertices of $P$ if $\underline{v}$ and $\overline{v}$ are finite. And the optimal solutions of the dual problems (\ref{primaljihe2}) and (\ref{primaljihe1dual}) remain unchange.

(2) If $\Psi(v)=[\underline{u},\overline{u}]$ is an interval, where $\underline{u}\ne\overline{u}$, then all dual optimal solutions $\hat{w}^*(u)$ of the problem (\ref{primaljihe1dual}) for $u\in [\underline{u},\overline{u}]$ forms an edge of the dual feasible region $D$, in which $\hat{w}^*(\underline{u})$ and $\hat{w}^*(\overline{u})$ are two adjoint vertices of $D$ if $\underline{u}$ and $\overline{u}$ are finite. And the optimal solutions of the primal problems (\ref{primaljihe1}) and (\ref{primaljihe2dual}) remain unchange. \upshape

\lem \label{AMyiny4} Let $r=1$.

(1) If $\Phi(u)=[\underline{v},\overline{v}]$ is an interval, where $\underline{v}\ne\overline{v}$, then for any $v\in (\underline{v},\overline{v})$, $\Psi(v)=u$.

(1) If $\Psi(v)=[\underline{u},\overline{u}]$ is an interval, where $\underline{u}\ne\overline{u}$, then for any $u\in (\underline{u},\overline{u})$, $\Phi(u)=v$. \upshape

Lemma \ref{AMyiny3} provides us with a pivot path for the solution of the problem (\ref{primaljihe10p}). From Lemma \ref{AMyiny3}, a projection interval of $\varTheta_P$ is the union of finite open invariancy intervals and transition points, see also Adler and Monterio \cite{AM92}. That is, there is a finite set of transition points $v_1<v_2<\dotsc<v_k$ such that
\[\varTheta_P=\bigcup\limits_{i=1}^{k-1}[v_i,v_{i+1}], \] in which every open interval $(v_i,v_{i+1})$ is an invariancy interval, and $v_1=-\infty$ and/or $v_k=+\infty$ are allowed. Hence every optimal solution $\hat{x}^*(v_i)$ of the problem (\ref{primaljihe2dual}) is a vertex of $P$, and every linear segment $[\hat{x}^*(v_i), \hat{x}^*(v_{i+1})]$ is an edge of $P$. Since $0\in \varTheta_P$, i.e., the problem (\ref{primaljihe2dual}) has an optimal solution $\hat{x}^*(0)$, the following path \[\hat{x}^*(v_1), \hat{x}^*(v_{2}),\cdots,\hat{x}^*(v_k) \] covers a pivot path for the solution of the problem (\ref{primaljihe10p}).

On the other hand, by Lemma \ref{AMyiny4}, for the dual projection interval $\varTheta_D$, there is a finite set of transition points $u_1<u_2<\dotsc<u_{k-1}$ such that $\Psi((v_i,v_{i+1}))=u_i$ for $i=1,2,\cdots,k-1$. And every $\Psi(v_i)$ is a closed interval whose interiors forms an open invariancy interval of $\varTheta_D$ if $v_i$ is finite. In particular, if $v_1$ is finite, then $\Psi(v_1)=[u_{k-1},+\infty)$; and if $v_k$ is finite, then $\Psi(v_k)=(-\infty,u_1]$. All in all, if $\varTheta_P$ contains $k-1$ open invariancy intervals, then $\varTheta_D$ contains $k-2$ to $k$ open invariancy intervals. Hence every optimal solution $\hat{w}^*(u_j)$ of the problem (\ref{primaljihe1dual}) is a vertex of $D$, and every linear segment $[\hat{w}^*(u_j), \hat{w}^*(u_{j+1})]$ is an edge of $D$. Similarly, we can obtain a dual pivot path by
 \[\hat{w}^*(u_1), \hat{w}^*(u_{2}),\cdots,\hat{w}^*(u_{k-1}) \] since $0\in \varTheta_D$.

\lem \label{finite2} If $r=1$ and $n-m-r=1$, then $\varTheta_P$ contains $n$ transition points and $n-1$ invariancy intervals at most. \upshape

\prof Let $B=(b_{ij})\in\mathbb{R}^{n\times l}$ be defined as in Corollary \ref{xinzd1}. If $n-m-r=1$, then $l=1$ such that the parametric LP (\ref{primaljihe2}) is easy to be solved by the use of the dual simplex algorithm since its dual has only one constraint. Hence the dual feasible basic variable $y_j$ can be chosen by the ratio test
\[ j\in \mathcal{J}=\left\{j\left| \frac{Bc}{b_{i1}}>0, b_{i1}\ne 0\right. \right\}. \] Clearly, the number of the entries of the index set $\mathcal{J}$ is less than or equal to $n$. Then by Lemma \ref{AMyiny3}, we conclude the final result. $\square$

The main result of this paper is as follows.

\thm\label{mainlp1} There is a pivot rule such that the simplex algorithm solves the problem (\ref{primaljihe10p}) at most in $n$ steps. \upshape

\prof For the given problem (\ref{primaljihe10p}), consider its single-parameter perturbation of the objective. Without loss of generality, we can assume that  for the primal problem (\ref{primaljihe10c}), $M^TS$ with rank 1 denotes only one perturbation direction, where $S$ is a $r\times r$ symmetric projection matrix. Meanwhile, we can also choose $r=n-m-1$ such that the dual problem (\ref{primaljihe10cd}) contains only one dual constraint (see also the proof of Lemma \ref{finite2}). If $r=0$, then the result is trivial; otherwise, we can choose $S$ with rank 1. By Lemma \ref{finite2} for the projection problem (\ref{primaljihe10cd}), the set $S\varTheta_P$ contains $n$ transition points and $n-1$ invariancy intervals at most. Then by Theorem \ref{mfzjia1}, the set $S\varTheta_D$ contains $n+1$ transition points and $n$ invariancy intervals at most. Finally, by Lemma \ref{AMyiny3}, the desired result is proved. $\square$

\end{document}